\documentclass[11pt]{article}

\usepackage{amscd,amsmath, amssymb, fancyhdr, epsfig,url}

\numberwithin{equation}{section}

\newcommand{\version}{version 4.1,\ \   March 18, 2015}

\def\eqref#1{(\ref{#1})}
\newcommand{\goth}{\mathfrak}

\newcommand{\arrow}{{\:\longrightarrow\:}}
\newcommand{\Z}{{\Bbb Z}}
\newcommand{\C}{{\Bbb C}}

\newcommand{\R}{{\Bbb R}}
\newcommand{\Q}{{\Bbb Q}}

\def\1{\sqrt{-1}\:}
\newcommand{\restrict}[1]{{\left|_{{\phantom{|}\!\!}_{#1}}\right.}}
\newcommand{\cntrct}                
{\hspace{2pt}\raisebox{1pt}{\text{$\lrcorner$}}\hspace{2pt}}

\makeatletter
\def\x@arrow{\DOTSB\Relbar}
\def\xlongequalsignfill@{\arrowfill@\x@arrow\Relbar\x@arrow}
\newcommand{\xlongequal}[2][]{%
        \ext@arrow 0099\xlongequalsignfill@{#1}{#2}}
\def\xlongrightarrowfill@{\arrowfill@\relbar\relbar\longrightarrow}
\newcommand{\xlongrightarrow}[2][]{%
        \ext@arrow 0099\xlongrightarrowfill@{#1}{#2}}
\makeatother

\renewcommand{\bar}{\overline}
\renewcommand{\phi}{\varphi}
\renewcommand{\epsilon}{\varepsilon}
\renewcommand{\geq}{\geqslant}
\renewcommand{\leq}{\leqslant}

\newcommand{\Gr}{\operatorname{Gr}}
\newcommand{\Hyp}{\operatorname{Hyp}}
\newcommand{\St}{\operatorname{\sf St}}
\newcommand{\End}{\operatorname{End}}

\newcommand{\Id}{\operatorname{Id}}

\newcommand{\Hom}{\operatorname{Hom}}

\newcommand{\Pic}{\operatorname{Pic}}
\newcommand{\Pos}{\operatorname{Pos}}
\newcommand{\Kah}{\operatorname{Kah}}

\newcommand{\Sym}{\operatorname{Sym}}

\newcommand{\Diff}{\operatorname{Diff}}

\newcommand{\rk}{\operatorname{rk}}

\newcommand{\Tw}{\operatorname{Tw}}

\renewcommand{\Re}{\operatorname{Re}}
\renewcommand{\Im}{\operatorname{Im}}

\newcommand{\Teich}{\operatorname{\sf Teich}}
\newcommand{\Comp}{\operatorname{\sf Comp}}

\newcommand{\Per}{\operatorname{\sf Per}}
\newcommand{\Perspace}{\operatorname{{\Bbb P}\sf er}}

\newcounter{Mycounter}[section]
\newcounter{lemma}[section]
\renewcommand{\thelemma}{{Lemma \thesection.\arabic{lemma}}}
\newcommand{\lemma}{%
    \setcounter{lemma}{\value{Mycounter}}
    \refstepcounter{lemma}
    \stepcounter{Mycounter}
    {\noindent \bf \thelemma:\ }}

\newcounter{claim}[section]
\renewcommand{\theclaim}{{Claim \thesection.\arabic{claim}}}
\newcommand{\claim}{%
    \setcounter{claim}{\value{Mycounter}}
    \refstepcounter{claim}
    \stepcounter{Mycounter}
    {\noindent \bf \theclaim:\ }}

\newcounter{sublemma}[section]
\newcounter{corollary}[section]
\renewcommand{\thecorollary}{{Corollary \thesection.\arabic{corollary}}}
\newcommand{\corollary}{%
    \setcounter{corollary}{\value{Mycounter}}
    \refstepcounter{corollary}
    \stepcounter{Mycounter}
    {\noindent \bf \thecorollary:\ }}

\newcounter{theorem}[section]
\renewcommand{\thetheorem}{{Theorem \thesection.\arabic{theorem}}}
\newcommand{\theorem}{%
    \setcounter{theorem}{\value{Mycounter}}
    \refstepcounter{theorem}
    \stepcounter{Mycounter}
    {\noindent \bf \thetheorem:\ }}

\newcounter{conjecture}[section]
\renewcommand{\theconjecture}{{Conjecture \thesection.\arabic{conjecture}}}
\newcommand{\conjecture}{%
    \setcounter{conjecture}{\value{Mycounter}}
    \refstepcounter{conjecture}
    \stepcounter{Mycounter}
    {\noindent \bf \theconjecture:\ }}

\newcounter{proposition}[section]
\renewcommand{\theproposition}
      {{Proposition \thesection.\arabic{proposition}}}
\newcommand{\proposition}{%
    \setcounter{proposition}{\value{Mycounter}}
    \refstepcounter{proposition}
    \stepcounter{Mycounter}
    {\noindent \bf \theproposition:\ }}

\newcounter{definition}[section]
\renewcommand{\thedefinition}
      {{Definition~\thesection.\arabic{definition}}}
\newcommand{\definition}{%
    \setcounter{definition}{\value{Mycounter}}
    \refstepcounter{definition}
    \stepcounter{Mycounter}
    {\noindent \bf \thedefinition:\ }}

\newcounter{example}[section]
\renewcommand{\theexample}{{Example \thesection.\arabic{example}}}
\newcommand{\example}{%
    \setcounter{example}{\value{Mycounter}}
    \refstepcounter{example}
    \stepcounter{Mycounter}
    {\noindent \bf \theexample:\ }}

\newcounter{remark}[section]
\renewcommand{\theremark}{{Remark \thesection.\arabic{remark}}}
\newcommand{\remark}{%
    \setcounter{remark}{\value{Mycounter}}
    \refstepcounter{remark}
    \stepcounter{Mycounter}
    {\noindent \bf \theremark:\ }}

\newcounter{problem}[section]
\newcounter{question}[section]
\makeatletter

\@addtoreset{equation}{section} \@addtoreset{footnote}{section} \makeatother

\def\blacksquare{\hbox{\vrule width 5pt height 5pt depth 0pt}}
\def\endproof{\blacksquare}

\begin{document}
\begin{center}
{\LARGE\bf
Ergodic complex structures on hyperk\"ahler manifolds\\[4mm]
}

Misha Verbitsky\footnote{Partially supported by 
RSCF grant 14-21-00053 within AG Laboratory NRU-HSE. }

\end{center}

{\small \hspace{0.10\linewidth}
\begin{minipage}[t]{0.85\linewidth}
{\bf Abstract} \\
Let $M$ be a compact complex manifold. The 
corresponding Teichm\"uller space $\Teich$ is a space
of all complex structures on $M$ up to the action
of the group $\Diff_0(M)$ if isotopies. The mapping
class group $\Gamma:=\Diff(M)/\Diff_0(M)$ acts
on $\Teich$ in a natural way. An {\bf ergodic complex
structure} is the one with a $\Gamma$-orbit dense 
in $\Teich$. Let $M$ be a complex torus of complex dimension
$\geq 2$ or a hyperk\"ahler manifold with $b_2>3$.
We prove that $M$ is ergodic, unless $M$ has maximal 
Picard rank (there is a countable number of such $M$).
This is used to show that all hyperk\"ahler manifolds
are Kobayashi non-hyperbolic.
\end{minipage}
}

\tableofcontents


\section{Introduction}


\subsection{Complex geometry and ergodic theory}

For an introduction to Teichm\"uller theory
and global Torelli theorem, please see
Subsection \ref{_Teich_intro_Subsection_};
the basic notions of hyperk\"ahler geometry
are recalled in Section \ref{_hk_Section_}. 
Here we assume that a reader knows the basic
definitions.

The Teichm\"uller space is defined as a space
of complex structures up to isotopies:
$\Teich:=\Comp/\Diff_0$. 
The mapping class group (also known as 
``a group of diffeotopies'') is the group
$\Gamma:= \Diff/\Diff_0$ of connected
components of the diffeomorphism group.
Clearly, $\Gamma$ acts on the Teichm\"uller
space in a natural way.

It turns out that in some important geometric situations
(for the hyperk\"ahler manifolds with $b_2>3$ and the complex tori of dimension
$\geq 2$) the mapping group action on $\Teich$ is ergodic.
This is surprising, at least to the author of 
the present paper, because in this case ``the moduli
space'' $\Teich/\Gamma$ of these geometric objects is extremely
pathological. In fact this quotient is so much non-Hausdorff that
all non-empty open subsets of $\Teich/\Gamma$  intersect
(\ref{_open_subse_intersect_Remark_}).

Complex structures with dense $\Gamma$-orbits
are called {\bf ergodic} (see \ref{_ergo_co_Definition_}).
From the description of the moduli in terms of 
homogeneous spaces and Moore's theorem on ergodic actions
it follows that the set
of non-ergodic complex structures on hyperk\"ahler 
manifolds with $b_2>3$ and complex tori of dimension
$\geq 2$ has measure 0
(\ref{_non_ergo_measure_0_Theorem_}). 
Applying Ratner theory, we prove that the set
of non-ergodic complex structures is in fact 
countable: a complex structure is non-ergodic
if and only if its Picard rank is maximal 
(\ref{_maxima_rank_non-ergo_Corollary_}). 

The density of particular
families of hyperk\"ahler manifolds in $\Teich/\Gamma$
was used many times since early 1970-ies.  
Pyatetsky-Shapiro and Shafarevich
used density of the family of Kummer surfaces in the
moduli of K3 surfaces to prove the local Torelli theorem
(\cite{_Piatetski_Shapiro_Shafarevich_}). 
This theorem was generalized to a general hyperk\"ahler
manifold $M$ with $b_2 \geq 5$ in 
\cite{_Ananin_Verbitsky_}. Here it was proven 
that any divisorial family defined by an
integer class in $H^2(M, \Z)$ is dense in
$\Teich/\Gamma$. In \cite{_Kamenova_V:fibrations_}
this approach was used further to study the 
Lagrangian  fibrations on hyperk\"ahler manifold.
Using the density argument and existence
of Lagrangian fibrations it was proven that 
all known hyperk\"ahler manifolds
are non-hyperbolic.\footnote{In the present paper
we generalize this further to all
hyperk\"ahler manifolds with $b_2 > 3$.}
In \cite{_Markman_Mehrota_},
Markman and Mehrota show that the space of 
Hilbert space of K3 surfaces is dense
in the corresponding deformation space,
and prove a similar result about the
generalized Kummer varieties.

Existence of ergodic complex structures leads to
some interesting
results about various complex-analytic 
quantities, such as the Kobayashi pseudometric.
As a model situation, consider a function 
$\phi$ on the set of equivalence classes of complex manifolds 
which is continuous
on deformations. Since the ergodic orbit
$\Gamma\cdot I$ is dense in the Teichm\"uller space,
and $\phi$ is constant on $\Gamma\cdot I$,
this implies that $\phi$ is constant.

In practice, such an application is hard to come
by, because functions which continuously depend
on the complex structure are sowewhat rare. However,
there are many semicontinuous functions, and
a semicontinuous function has to be constant
on ergodic complex structures. Indeed, let $I, J\in \Teich$
be two ergodic complex structures, that is,
complex structures with with dense $\Gamma$-orbits,
and $\phi:\; \Teich\arrow \R$ a semicontinuous
(say, upper semicontinuous) $\Gamma$-invariant function.
Since $I$ is a limiting point of a dense set 
$\Gamma\cdot J$, semicontinuity implies $\phi(I)\geq \phi(J)$.
By the same reason, $\phi(J)\geq \phi(I)$,
hence $\phi$ is constant on the set of all
ergodic complex structures. 

This observation can be applied to several
questions of complex hyperbolicity 
(Subsection \ref{_kobaya_intro_Subsection_}).

\subsection{Teichm\"uller spaces and hyperk\"ahler
  geometry}
\label{_Teich_intro_Subsection_}

We recapitulate briefly the definition of a Teichm\"uller space
of the hyper\-k\"ahler manifolds, following \cite{_V:Torelli_}.

\hfill

\definition
Let $M$ be a compact complex manifold, and 
$\Diff_0(M)$ a connected component of its diffeomorphism group
({\bf the group of isotopies}). Denote by $\Comp$
the space of complex structures on $M$, equipped with
a structure of Fr\'echet manifold, and let
$\Teich:=\Comp/\Diff_0(M)$ be its quotient, equipped
with the quotient topology. We call 
it {\bf the Teichm\"uller space.}

\hfill

\remark
In many important cases, such as
for Calabi-Yau manifolds (\cite{_Catanese:moduli_}), 
$\Teich$ is a finite-dimensional
complex space; usually it is non-Hausdorff.

\hfill

\definition
Let $\Diff(M)$ be the group of orientable diffeomorphisms of
a complex manifold. The quotient $\Comp/\Diff=\Teich/\Gamma$ 
is called {\bf the moduli space} of complex structures on $M$.
Typically, it is very non-Hausdorff. The set
$\Comp/\Diff$ corresponds bijectively to the set of isomorphism
classes of complex structures.

\hfill

\definition
A {\bf hyperk\"ahler structure} on a manifold $M$
is a Riemannian structure $g$ and a triple of complex
structures $I,J,K$, satisfying quaternionic relations
$I\circ J = - J \circ I =K$, such that $g$ is K\"ahler
for $I,J,K$.

\hfill

\remark
One could define a hyperk\"ahler structure in 
terms of the complex geometry of its 
twistor space (\ref{_twistor_Definition_}). 
This was discovered in \cite{_HKLR_}; 
see \cite{_Verbitsky:hypercomple_}
for a few historical remarks and further development
of this approach.

\hfill

\remark
 A hyperk\"ahler manifold is holomorphically
symplectic: $\omega_J+\1 \omega_K$ is a holomorphic
symplectic form on $(M,I)$, which is easily seen
using a simple linear-algebraic calculation
(\cite{_Besse:Einst_Manifo_}).

\hfill

\theorem \label{_Calabi_Yau_Theorem_}
(Calabi-Yau; see \cite{_Yau:Calabi-Yau_,_Beauville_,_Besse:Einst_Manifo_}) \\
A compact, K\"ahler, holomorphically symplectic manifold
admits a unique hyperk\"ahler metric in any K\"ahler class.
\endproof

\hfill

\remark
The term ``hyperk\"ahler manifold'' can mean
many different things. In the literature, it
denotes either a manifold equipped with a hyperk\"ahler
structure, or a complex manifold admitting a hyperk\"ahler
structure, or a Riemannian manifold with holonomy
in $Sp(n)$. In the present paper, we shall understand
by a ``hyperk\"ahler manifold'' a compact, complex
manifold admitting a K\"ahler structure and 
a holomorphically symplectic structure.
Such a complex structure is called
a complex structure of hyperk\"ahler type. We also
assume tacitly that all hyperk\"ahler manifolds
are {\em of maximal holonomy} (or ``irreducibly holomorphically
symplectic''), in the sense of \ref{_hk_simple_Definition_} 
below.

\hfill

\definition\label{_hk_simple_Definition_}
A hyperk\"ahler manifold $M$ is called
{\bf of maximal holonomy} if $\pi_1(M)=0$, $H^{2,0}(M)=\C$.
In the literature, such manifolds are often
called {\bf irreducibly holomorphic symplectic},
or {\bf irreducibly symplectic varieties}.

\hfill

This definition is motivated by the following theorem
of Bogomolov. 

\hfill

\theorem (\cite{_Bogomolov:decompo_,_Beauville_})
Any hyperk\"ahler manifold admits a finite covering
which is a product of a torus and several 
maximal holonomy hyperk\"ahler manifolds.
\endproof

\hfill

\remark\label{_Teich_hk_type_Remark_}
Further on, all hyperk\"ahler manifolds
are assumed to be of maximal holonomy, $\Comp$ is the
space of all complex structures of
hyperk\"ahler type on $M$, and $\Teich$
its quotient by $\Diff_0(M)$.

\subsection{Ergodic complex structures}

The main object of this paper is the following notion.

\hfill

\definition\label{_ergo_co_Definition_}
Let $M$ be a complex manifold, $\Teich$ its Teichm\"uller space, and
$I\in \Teich$ a point. Consider the set $Z_I\subset \Teich$
of all $I'\in \Teich$ such that $(M,I)$ is biholomorphic to $(M,I')$
(clearly, $Z_I=\Gamma \cdot I$, where $\Gamma=\Diff(M)/\Diff_0(M)$ 
is a mapping group acting on
$\Teich$). A complex structure is called {\bf ergodic}
if the corresponding orbit $Z_I$ is dense in $\Teich$.

\hfill

\remark
The origins of this term are explained in 
Subsection \ref{_ergo_Subsection_}.
It is well known that almost all orbits of
an ergodic action are dense (\ref{_non_dense_m_zero_Claim_}). 
However, for a hyperk\"ahler manifold with $b_2 >3$ 
and a complex torus of dimension $\geq 2$, the mapping group action
on $\Teich$ is ergodic (\ref{_non_ergo_measure_0_Theorem_},
\ref{_Teich_ergo_Remark_}). Notice that the definition
of ergodic action does not require one 
to fix a particular measure on the space
(\ref{_ergodic_on_mani_Remark_}).

\hfill

In many situations, the mapping class group action on the
Teichm\"uller space is ergodic. This implies
that the non-ergodic complex structures form 
a set of measure 0.

\hfill

\theorem\label{_non-ergo_pts_Theorem_}
Let $M$ be a maximal holonomy hyperk\"ahler manifold
or a compact complex torus of dimension $\geq 2$.
Then the set $\Teich_{ne}$ of non-ergodic points
has measure 0 in the corresponding Teichm\"uller 
space $\Teich$.

\hfill

{\bf Proof:} See \ref{_non_ergo_measure_0_Theorem_}. \endproof

\hfill

\remark
The notion of a measure zero subset of a manifold is 
independent from a choice of a smooth measure. Therefore,
to state \ref{_non-ergo_pts_Theorem_}, it is not necessary 
to fix a particular measure on $\Teich$.

\hfill

This result follows from the ergodicity of the
mapping class group action on $\Teich$, which
follows from the global Torelli theorem and the
ergodicity of an arithmetic action on homogeneous
spaces due to Calvin Moore (1966). It is not very explicit,
and for a considerable period of time, no explicit
examples of ergodic complex structures were known.
This problem was solved by an application of a powerful
theorem of Marina Ratner (\ref{_ratner_orbit_clo_Theorem_}).

\hfill

\theorem
Let $M$ be a maximal holonomy hyperk\"ahler manifold
or a compact complex torus of dimension $\geq 2$,
and $I$ a complex structure on $M$. Then $I$ is
non-ergodic iff the Neron-Severi lattice of $(M,I)$  
has maximal possible rank. This means that
$\rk NS(M,I)=b_2(M)-2$ for $M$ hyperk\"ahler, and
$\rk NS(M,I)=(\dim_\C M)^2$ for $M$ a torus.

\hfill

{\bf Proof:} See \ref{_maxima_rank_non-ergo_Corollary_}. \endproof

\subsection{Kobayashi pseudometric on hyperk\"ahler manifolds}
\label{_kobaya_intro_Subsection_}

The ergodic properties of the mapping group action
have many applications to Kobayashi hyperbolicity.
For the definition of Kobayashi pseudometric, 
basic properties and the further reference, please see Section 
\ref{_Kobayashi_Section_}. For the purposes of the
present paper, Kobayashi pseudometric is important
because it is a complex-analytic invariant which
is upper semicontinuous as a function of a 
complex structure (\cite{_Kobayashi:1970_,_Voisin:kobayashi_}).
This suggests the following conjecture.

\hfill

\conjecture\label{_Kobaya_iso_Conjecture_}
Let $I, J\in \Teich$ be ergodic complex structures
on $M$, and $d_I, d_J$ the corresponding
Kobayashi pseudometrics. Then $(M, d_I)$
is isometric to $(M,d_J)$.

\hfill

\remark
Since the $\Gamma$-orbit of $I$ is dense in $\Teich$,
for any $K\in \Teich$, $K$ can be obtained as a limit
of $\nu_i^* I$ (\ref{_ergo_limits_Remark_}), and 
one has $d_K\geq  d_I$ by semicontinuity of $K$.
In principle, this should give  $d_J \geq d_I \geq d_J$, because
both $I$ and $J$ are ergodic.
To make this heuristic argument rigorous, one should make
the dependency of $d_I$, $d_J$ on diffeomorpisms 
$\nu_i$ explicit.

\hfill

Kobayashi hyperbolic manifolds are those with 
non-degenerate Ko\-ba\-ya\-shi pseudometric. The set of 
Kobayashi hyperbolic complex structures is 
open in holomorphic families. Moreover, given 
a holomorphic family of Kobayashi hyperbolic manifolds, the Kobayashi
pseudometric is continuous in this family 
(\cite{_Voisin:kobayashi_}). In the presence
of an ergodic complex structure, the argument used
in the sketch of \ref{_Kobaya_iso_Conjecture_}
would imply that all hyperbolic complex structures
on $M$ are isometric, and the mapping class group acts
by isometries. However, the isometry group of a compact
metric space is always compact, and the image of the
mapping class group in cohomology is usually non-compact.
This can be used to prove non-hyperbolicity
for manifolds admitting ergodic complex structures.

However, for hyperk\"ahler manifolds, there exists
a very simple and direct argument proving non-hyperbolicity.

The non-hyperbolicity of hyperk\"ahler (and, more
generally, Calabi-Yau) manifolds was a subject
of long research, but until recently the only 
general result was a theorem by F. Campana
proven in \cite{_Campana:twistor_hy_}. 

A {\bf twistor space} of a hyperk\"ahler manifold
(\ref{_twistor_Definition_}) is a total space of a fibration
obtained from a hyperk\"ahler rotation of a complex
structure. Campana observed that the space of
rational curves on $\Tw(M)$ transversal to the fibers 
is never compact (in fact it is holomorphically convex, as shown in
\cite{_NHYM_}; see also \cite{_Verbitsky:special_twistor_}
and \cite{_DD_Mourougane_}). 
Then a limit of a sequence of rational curves in $\Tw(M)$
would contain an entire curve in one of the twistor
fibers. This implies the Campana non-hyperbolicity theorem
(\ref{_Campana_non_hy_Theorem_}): at least one of the
fibers of the twistor family is not Kobayashi hyperbolic.

Another approach to non-hyperbolicity was used in
\cite{_Kamenova_V:fibrations_}. In this paper it was shown
that for all known examples of hyperk\"ahler manifolds,
the manifolds admitting a holomorphic Lagrangian fibration
are dense in the moduli space. Such manifolds 
contain entire curves, hence they are non-hyperbolic.
However, the set of non-hyperbolic complex structures
is closed in the relevant deformation space by Brody's lemma
(\ref{_limi_non-hy_Theorem_}). This is how
non-hyperbolicity was proven in
\cite{_Kamenova_V:fibrations_}.

In the present paper we go by a different route, using
ergodic methods and the Campana's non-hyperbolicity
result. The explicit description of the set of non-ergodic
complex structures (\ref{_maxima_rank_non-ergo_Corollary_})
allows one to find a twistor family with all fibers
ergodic (\ref{_twi_ergo_Claim_}).
By Campana's theorem, one of these fibers
is non-hyperbolic. This gives a non-hyperbolic
ergodic complex structure $I\in \Teich$.
Then all points of the set $\Gamma\cdot I\subset \Teich$
are also non-hyperbolic. However, this set is dense,
hence its closure $\overline{\Gamma\cdot I}$ is the whole of $\Teich$.
Finally, we observe that the set of non-hyperbolic
complex structures is closed in $\Teich$, and therefore
it contains $\overline{\Gamma\cdot I}=\Teich$.


\section{Hyperk\"ahler manifolds}
\label{_hk_Section_}


In this section, we state  the global Torelli theorem for
hyperk\"ahler manifolds, following \cite{_V:Torelli_}.

\subsection{Bogomolov-Beauville-Fujiki form}

The Bogomolov-Beauville-Fujiki form was
defined in \cite{_Bogomolov:defo_} and 
\cite{_Beauville_},
but it is easiest to describe it using the
Fujiki formula, proven in \cite{_Fujiki:HK_}.

\hfill

\theorem\label{_Fujiki_Theorem_}
(Fujiki)
Let $M$ be a maximal holonomy hyperk\"ahler manifold,
$\eta\in H^2(M)$, and $n=\frac 1 2 \dim M$. 
Then $\int_M \eta^{2n}=c q(\eta,\eta)^n$,
where $q$ is a primitive integral quadratic form on $H^2(M,\Z)$,
and $c>0$ a rational number. \endproof

\hfill

\remark 
Fujiki formula (\ref{_Fujiki_Theorem_}) 
determines the form $q$ uniquely up to a sign.
For odd $n$, the sign is unambiguously determined as well.
For even $n$, one needs the following explicit
formula, which is due to Bogomolov and Beauville.
\begin{equation}\label{_BBF_expli_Equation_}
\begin{aligned}  \lambda q(\eta,\eta) &=
   \int_X \eta\wedge\eta  \wedge \Omega^{n-1}
   \wedge \bar \Omega^{n-1} -\\
 &-\frac {n-1}{n}\left(\int_X \eta \wedge \Omega^{n-1}\wedge \bar
   \Omega^{n}\right) \left(\int_X \eta \wedge \Omega^{n}\wedge \bar \Omega^{n-1}\right)
\end{aligned}
\end{equation}
where $\Omega$ is the holomorphic symplectic form, and 
$\lambda>0$.

\hfill

\definition
Let $q\in \Sym^2(H^2(M,\Z)^*)$ be the integral form
defined by \ref{_Fujiki_Theorem_} and
\eqref{_BBF_expli_Equation_}.
This form is called {\bf the Bogomolov-Beauville-Fujiki
  form}.

\subsection{Mapping class group}

\definition
 Let $\Diff(M)$ be the group of oriented
diffeomorphisms of $M$, and $\Diff_0(M)$ the
group of isotopies, that is, the
connected component of $\Diff(M)$. 
We call $\Gamma:=\Diff(M)/\Diff_0(M)$ {\bf the
mapping class group} of $M$.

\hfill

For K\"ahler manifolds of dimension $\geq 3$, 
the mapping class group can be computed using the
following theorem of Sullivan.

\hfill

\theorem \label{_Sullivan_mapping_class_Theorem_}
(D. Sullivan; \cite{_Sullivan:infinite_}) \\
Let $M$ be a compact, simply connected 
K\"ahler manifold, $\dim_\C M\geq 3$. Denote by $\Gamma_0$ the group
of automorphisms of the algebra $H^*(M, \Z)$
preserving the Pontryagin classes $p_i(M)$. 
Then  the natural map 
$\Diff(M)/\Diff_0\arrow \Gamma_0$ has finite kernel,
and its image has finite index in $\Gamma_0$.
\endproof

\hfill

\definition
Groups $G, G'$ are {\bf commensurable}
if there exist subgroups $G_1\subset G, G_1' \subset G_1$
of finite index, and finite normal subgroups $G_2\subset G_1$
and $G_2' \subset G_1'$ such that $G_1/G_2$ is isomorphic
to $G_1'/G_2'$. {\bf An arithmetic group}
is a group which is commensurable to an integer lattice
in a rational Lie group. 

\hfill

\remark 
Sullivan's theorem claims that the mapping
class group of any K\"ahler manifold is
arithmetic.\footnote{In fact, Sullivan proved the
arithmeticity of the mapping class group
for any compact smooth manifold of dimension $\geq 5$.}

\hfill

Using the results of \cite{_Verbitsky:coho_announce_},
the group of automorphisms of the algebra $H^*(M, \Z)$
can be determined explicitly, up to commensurability.
This gives the following theorem, proven in 
\cite{_V:Torelli_}.

\hfill

\theorem\label{_mapping_class_Theorem_}
Let $M$ be a maximal holonomy hyperk\"ahler manifold,
and $\Gamma_0$ the group of automorphisms of $H^*(M,\Z)$
preserving Pontryagin classes 
(\ref{_Sullivan_mapping_class_Theorem_}).
Consider the restriction map $\Gamma_0 \stackrel \psi\arrow GL(H^2(M, \Z))$.
Then $\psi$ has finite kernel, its image lies in
the orthogonal group $O(H^2(M, \Z), q)$, and $\psi(\Gamma_0)$
has finite index in $O(H^2(M, \Z), q)$.
\endproof

\subsection{Global Torelli theorem}

\remark
Let $M$ be a hyperk\"ahler manifold (as usually, we assume
$M$ to be of maximal holonomy). Recall that in this situation
$\Teich$ was defined as a set of all complex
structures of hyperk\"ahler type on $M$ (\ref{_Teich_hk_type_Remark_}).
For any $J$ in the same connected components of $\Teich$,
$(M,J)$ is also a maximal holonomy hyperk\"ahler manifold, 
because the Hodge numbers are constant in families.
Therefore, $H^{2,0}(M,J)$ is one-dimensional. 

\hfill

\definition
 Let 
\[ \Per:\; \Teich \arrow {\Bbb P}H^2(M, \C)
\]
map $J$ to the line $H^{2,0}(M,J)\in {\Bbb P}H^2(M, \C)$.
The map $\Per$ is 
called {\bf the period map}.

\hfill

\remark
The period map $\Per$ maps $\Teich$ into an open subset of a 
quadric, defined by
\[
\Perspace:= \{l\in {\Bbb P}H^2(M, \C)\ \ | \ \  q(l,l)=0, q(l, \bar l) >0\}.
\]
It is called {\bf the period space} of $M$.
Indeed, any holomorphic symplectic form $l$
satisfies the relations $q(l,l)=0, q(l, \bar l) >0$,
as follows from \eqref{_BBF_expli_Equation_}.

\hfill

\proposition\label{_period_Grassmann_Proposition_}
The period space $\Perspace$
is identified with the quotient
$SO(b_2-3,3)/SO(2) \times SO(b_2-3,1)$, which
is a Grassmannian $\Gr_{++}(H^2(M,\R))$
of positive oriented 2-planes in $H^2(M,\R)$.

\hfill

{\bf Proof:} This statement is well known, but we shall
sketch its proof to illustrate the constructions given
below.

{\bf Step 1:} Given $l\in {\Bbb P}H^2(M, \C)$, the space
generated by $\Im l, \Re l$ is 2-dimensional, because 
$q(l,l)=0, q(l, \bar l)>0$ implies that $l \cap H^2(M,\R)=0$.

{\bf Step 2:}  This 2-dimensional plane is 
positive, because 
 $q(\Re l, \Re l) = q(l+ \bar l, l+ \bar l) = 2 q(l, \bar l)>0$.

{\bf Step 3:} Conversely, for any 2-dimensional positive
plane  $V\in H^2(M,\R)$, 
the quadric $\{l\in V \otimes_\R \C\ \ | \ \ q(l,l)=0\}$
consists of two lines; a choice of a line is determined by orientation.
\endproof

\hfill

\definition
Let $M$ be a topological space. We say that $x, y \in M$
are {\bf non-separable} (denoted by $x\sim y$)
if for any open sets $V\ni x, U\ni y$, $U \cap V\neq \emptyset$.

\hfill

\theorem \label{_Huy_non_separa_Theorem_}
(Huybrechts; \cite{_Huybrechts:basic_}).
Two points $I,I'\in \Teich$ are non-separable if  
there exists a bimeromorphism $(M,I)\arrow (M,I')$.
\endproof

\hfill

\definition
The space $\Teich_b:= \Teich\!/\!\!\sim$ is called {\bf the
birational Teichm\"uller space} of $M$.

\hfill

\theorem 
(Global Torelli theorem; \cite{_V:Torelli_})
The period map 
$\Teich_b\stackrel \Per \arrow \Perspace$ is an isomorphism,
for each connected component of $\Teich_b$.
\endproof

\hfill

\definition
Let $M$ be a hyperk\"ahler manifold,
$\Teich_b$ its birational Teichm\"uller space,
and $\Gamma$ the mapping class group.
The quotient $\Teich_b/\Gamma$ is called
{\bf the birational moduli space} of $M$.
Its points are in bijective correspondence with the
complex structures of hyperk\"ahler type on $M$
up to a bimeromorphic equivalence.

\hfill

\remark
The word ``space'' in this context is misleading.
In fact, the quotient topology on $\Teich_b/\Gamma$ is extremely
non-Hausdorff, e.g. every two open sets would intersect
(\ref{_open_subse_intersect_Remark_}).

\hfill

The Global Torelli theorem can be stated
as a result about the birational moduli space.

\hfill

\theorem\label{_moduli_monodro_Theorem_}
(\cite[Theorem 7.2, Remark 7.4]{_V:Torelli_})
Let $(M,I)$ be a hyperk\"ahler manifold, and $W$ 
a connected component of its birational
moduli space. Then $W$ is isomorphic to ${\Perspace}/\Gamma_I$,
where ${\Perspace}=SO(b_2-3,3)/SO(2) \times SO(b_2-3,1)$
and $\Gamma_I$ is 
an arithmetic group in $O(H^2(M, \R), q)$, called {\bf the
monodromy group} of $(M,I)$.
\endproof

\hfill

\remark
The monodromy group of $(M,I)$ can be also described
as a subgroup of the group $O(H^2(M, \Z), q)$
generated by the monodromy transform maps for 
the Gauss-Manin local systems obtained from all
deformations of $(M,I)$ over a complex base
(\cite[Definition 7.1]{_V:Torelli_}). This is 
how this group was originally defined by Markman
(\cite{_Markman:constra_}, \cite{_Markman:survey_}).

\hfill

\remark 
A caution: usually ``the global Torelli theorem''
is understood as a theorem about Hodge structures.
For K3 surfaces, the Hodge structure on $H^2(M,\Z)$
determines the complex structure. 
For $\dim_\C M >2$, it is false.

\hfill

\remark
Further on, we shall freely identify $\Perspace$ and
$\Teich_b$.

\hfill

\remark\label{_fiber_Per_counta_Remark_}
By \cite[Proposition 5.14]{_Markman:survey_},
the fibers of the natural projection
$\Per:\; \Teich\arrow \Teich_b$ can be identified with 
a set of ``K\"ahler chambers'', which are 
open subsets of the space $H^{1,1}(M,I)$.
Therefore, each fiber is countable or finite.

\hfill

\remark\label{_fiber_Per_over_divi_Remark_}
By \cite[Remark 4.28]{_V:Torelli_},
outside of a countable union of complex divisors on $\Teich_b$,
the map $\Per:\; \Teich\arrow \Teich_b$ is bijective.

\hfill

\remark\label{_Teich_and_Per_measure_Remark_}
Further on, we would be interested in ergodic (that is,
measure-theoretic) properties of $\Teich$ and $\Teich_b$.
By \ref{_fiber_Per_over_divi_Remark_} and 
\ref{_fiber_Per_counta_Remark_}, the map $\Per$
is bijective outside of a measure 0 set.
Therefore, any ergodic theory result proven for $\Teich$
remains true for $\Teich_b$, and vice versa.

\section{Ergodic complex structures \\on 
hyperk\"ahler manifolds and tori}

\subsection{Ergodicity: basic definitions and results}
\label{_ergo_Subsection_}

\definition
Let $(M,\mu)$ be a space with measure,
and $G$ a group acting on $M$ preserving the sigma-algebra of measurable
subsets, and mapping measure zero sets to measure zero
sets. This action is {\bf ergodic} if all
$G$-invariant measurable subsets $M'\subset M$
satisfy $\mu(M')=0$ or $\mu(M\backslash M')=0$.

\hfill

\remark\label{_ergodic_on_mani_Remark_}
When one defines an ergodic action, it is usually assumed
that the action of $G$ preserves the measure. However,
this is not necessary. In fact, any manifold is equipped with
a sigma-algebra of Lebesgue measurable sets and, moreover, the notion of
a measure zero subset set is independent from a choice of a Lebesgue 
measure. This means that one can define ``ergodic action of a group''
on a manifold not specifying the measure.

\hfill

\claim\label{_non_dense_m_zero_Claim_}
Let $M$ be a manifold, $\mu$ a Lebesgue measure, and
$G$ a group acting on $(M,\mu)$ ergodically.  Then the 
set of points with non-dense orbits has measure 0.

\hfill

{\bf Proof:}
Consider a non-empty open subset $U\subset M$. 
Then $\mu(U)>0$, hence $M':=G\cdot U$ satisfies 
$\mu(M\backslash M')=0$. For any orbit $G\cdot x$
not intersecting $U$, $x\in M\backslash M'$.
Therefore the set of such points has measure 0.
\endproof

\hfill

\definition
Let $M$ be a complex manifold, $\Teich$ its Techm\"uller
space, and $\Gamma$ the mapping group acting on $\Teich$.
{\bf An ergodic complex structure} is a complex
structure with dense $\Gamma$-orbit.

\hfill

\remark\label{_ergo_limits_Remark_}
Let $(M,I)$ be a manifold with ergodic complex structure,
and $I'$ another complex structure.
Then there exists a sequence of diffeomorphisms
$\nu_i$ such that $\nu_i^*(I)$ converges to $I'$ in 
the usual (Fr\'echet) topology on the space of complex 
structure tensors. This property is clearly equivalent 
to ergodicity of $I$.

\hfill

Further on, we shall need the following result about
ergodicity of an arithmetic group action on a homogeneous
space. This result would be applied to a mapping class
group (which is arithmetic by \ref{_mapping_class_Theorem_}) and 
a period space, which is homogeneous 
(\ref{_period_Grassmann_Proposition_}).

\hfill

\definition 
Let $G$ be a Lie group, and $\Gamma\subset G$ a discrete
subgroup. Consider the pushforward of the Haar measure to
$G/\Gamma$. We say that $\Gamma$ {\bf has finite covolume}
if the Haar measure of $G/\Gamma$ is finite.
In this case $\Gamma$ is called {\bf a lattice subgroup}.

\hfill

\remark
Borel and Harish-Chandra proved that
an arithmetic subgroup of a reductive group $G$ over $\Q$
is a lattice whenever $G$ has no non-trivial characters
over $\Q$ (see e.g. \cite{_Vinberg_Gorbatsevich_Shvartsman_}). 
In particular, all arithmetic subgroups
of a semi-simple group defined over $\Q$ are lattices.

\hfill

\theorem \label{_Moore_Theorem_}
(Calvin C. Moore, \cite[Theorem 7]{_Moore:ergodi_})
Let $\Gamma$ be a lattice subgroup 
(such as an arithmetic subgroup) in a non-compact 
simple Lie group $G$ with finite center, and $H\subset G$ a 
non-compact Lie subgroup. Then the left action of $\Gamma$
on $G/H$ is ergodic. \endproof

\subsection{Ergodic action on the Teichm\"uller space for
hyperk\"ahler manifolds and tori}
\label{_ergo_complex_str_Subsection_}

\theorem\label{_non_ergo_measure_0_Theorem_}
 Let $\Perspace$ be a component of 
a birational Teichm\"uller space of a hyperk\"ahler
manifold $M$, $b_2(M)>3$,  and
$\Gamma_I$ its monodromy group acting on $\Perspace$. Consider the set
$Z\subset \Perspace$ of all  points with non-dense orbits.
Then the action of $\Gamma_I$ on $\Perspace$ is 
ergodic, and $Z$ has measure 0 in $\Perspace$.

\hfill

{\bf Proof. Step 1:} 
 Let $G=SO(b_2-3,3)$, $H=SO(2) \times SO(b_2-3,1)$,
and $\Gamma\subset G$ an arithmetic subgroup.
Then $\Gamma$-action on $G/H$ is ergodic, by Moore's theorem.

{\bf  Step 2:} The space $\Perspace$ is identified with
$G/H$ (\ref{_period_Grassmann_Proposition_}),
and the monodromy group is an arithmetic subgroup of $G$
by \ref{_mapping_class_Theorem_} and \ref{_moduli_monodro_Theorem_}.
Then $\Gamma_I$ acts on $\Perspace$ ergodically, and
the set of points with non-dense orbits has measure 0
(\ref{_non_dense_m_zero_Claim_}).
\endproof

\hfill

\remark\label{_Teich_ergo_Remark_}
As explained in \ref{_Teich_and_Per_measure_Remark_}, 
the space $\Teich_b=\Perspace$ is identified with $\Teich$
up to measure 0 subsets. Therefore, the set of non-ergodic complex structures
on a hyperk\"ahler manifolds has measure 0 in $\Teich$.

\hfill

A similar result is true for a compact torus.
Here the Teichm\"uller space is the space of complex
structure operators on $\R^{2n}$, identified with the quotient
$SL(2n,\R)/SL(n,\C)$, and the mapping class group is
$SL(2n, \Z)$ (\cite{_Catanese:moduli_}). For $n>1$,
the group $SL(n,\C)$ is non-compact.
Therefore, \ref{_Moore_Theorem_} can be applied, and 
we obtain the following statement.

\hfill

\theorem
Let $W:=SL(2n,\R)/SL(n,\C)$ be the Teichm\"uller
space of an $n$-dimensional compact torus, 
$n \geq 2$, equipped with an action of
the mapping class group $\Gamma=SL(2n, \Z)$.
Then the action of $\Gamma$ on $W$ is ergodic.
In particular, the set of non-ergodic complex
tori has measure 0 in the corresponding Teichm\"uller space.
\endproof

\hfill

\remark\label{_open_subse_intersect_Remark_}
Existence of erdogic complex structures means
that the quotient $\Teich\!/\Gamma$ (considered
with the quotient topology) is extremely
non-Hausdorff. Indeed, any two open sets in $\Teich$
contain points in a dense orbit $\Gamma\cdot I$,
hence their images in $\Teich\!/\Gamma$ intersect.
We obtain that any two open subsets in 
the moduli ``space'' $\Teich\!/\Gamma$ intersect.


\section{Ratner orbit closure theorem and ergodic complex structure}


\subsection{Lie groups generated by unipotents}

Here we state basic facts of Ratner theory.
We follow \cite{_Kleinbock_etc:Handbook_} 
and \cite{_Morris:Ratner_}.

\hfill

\definition
Let $G$ be a Lie group, and $g\in G$ any element.
We say that $g$ is {\bf unipotent} if $g=e^h$ for a
nilpotent element $h$ in its Lie algebra.
A group $G$ is {\bf generated by unipotents}
if $G$ is multiplicatively generated by unipotent elements.

\hfill

\theorem\label{_ratner_orbit_clo_Theorem_}
(Ratner orbit closure theorem)\\
Let $H\subset G$ be a Lie subroup generated by 
unipotents, and $\Gamma\subset G$ a lattice.
Then the closure of any $H$-orbit in $G/\Gamma$
is the orbit of a closed, connected subgroup $S\subset G$,
such that $S\cap \Gamma\subset S$ is a lattice in $S$.

\hfill

{\bf Proof:} \cite[1.1.15 (2)]{_Morris:Ratner_}. \endproof

\hfill

\remark
\ref{_ratner_orbit_clo_Theorem_} is 
true if $H=H_0\times H_1$,
where $H_0$ is generated by unipotents, and
$H_1$ is compact. Indeed, for each $x\in G/\Gamma$, one has
$\overline {H\cdot x}= H_1 \cdot \overline {H_0\cdot x}$.
The inclusion $\overline {H\cdot x}\supset 
H_1 \cdot \overline {H_0\cdot x}$ is obvious.
The converse inclusion would follow if we prove that
$H_1 \cdot \overline {H_0\cdot x}$ is closed. However,
the orbit of a closed set under a compact Lie group
is always closed.

\hfill

\example\label{_ergo_Example_}
Let $V$ be a real vector space with a non-degenerate
bilinear symmetric form of signature $(3,k)$, $k>0$
$G:=SO^+(V)$ a connected component of the 
isometry group, $H\subset G$ a subgroup
fixing a given positive 2-dimensional plane,
$H\cong SO^+(1,k)\times SO(2)$, and $\Gamma\subset G$ an arithmetic
lattice. Consider the quotient 
$\Perspace:=H\backslash G$. Then
\begin{description}
\item[(i)] A point $J\in \Perspace$ has closed $\Gamma$-orbit 
if and only if the orbit $H\cdot J$ in the quotient 
$G/\Gamma$ is closed.
\item[(ii)]   The closure of $H\cdot J$ 
in $G/\Gamma$ is the orbit of a closed connected
Lie group $S\supset H$:
\[ \overline{H\cdot J}= S\cdot J\subset \Perspace.\]
\end{description}

\hfill

For arithmetic groups
Ratner orbit closure theorem can be stated in a more
precise way, as follows.

\hfill

\theorem\label{_closure_arithm_Ratner_Theorem_}
Let $G$ be a real algebraic group defined over $\Q$
and with no non-trivial characters, $W\subset G$ a
subgroup generated by unipotents, and $\Gamma\subset G$
an arithmetic lattice. For a given $g\in G$,
let $H$ be the smallest real algebraic $\Q$-subgroup of
$G$ containing $g^{-1}Wg$. Then the closure of 
$Wg$ in $G/\Gamma$ is $Hg$. 

\hfill

{\bf Proof:} See \cite[Proposition 3.3.7]{_Kleinbock_etc:Handbook_}
or \cite[Proposition 3.2]{_Shah:uniformly_}.
\endproof

\subsection{Ratner theorem for Teichm\"uller spaces}

In Subsection \ref{_lie_subgrou_Subsection_}, we prove the following
two elementary theorems.

\hfill

\theorem\label{_subgro_SO(3,k)_Theorem_}
Let $G=SO^+(3,k)$, $k\geq 1$, and 
$H\cong SO^+(1,k)\times SO(2)\subset G$.
Then any closed connected
Lie subgroup $S\subset G$ containing $H$ coincides
with $G$ or with $H$.

{\bf Proof:} See \ref{_SO(p,q)_subgroup_Theorem_}. \endproof

\hfill

\theorem
Let $n\geq 2$, $G=SL(2n, \R)$, and 
$H\cong SL(n, \C) \subset G$.
Then any closed connected
Lie subgroup $S\subset G$ containing $H$ coincides
with $G$ or with $H$.

{\bf Proof:} See \ref{_complex_automo_maxi_Theorem_}. \endproof

\hfill

Now we can apply these theorems to characterize 
the ergodic and non-ergodic complex structures.

\hfill

\theorem\label{_perspace_dense_Theorem_}
Let $M$ be a hyperk\"ahler manifold, 
$\Perspace$ its period space, and 
$I\in \Perspace=Gr_{++}(H^2(M, \R))$ a point associated to
a positive 2-plane $V\subset H^2(M, \R)$. Then 
the $\Gamma$-orbit of $I$ is dense is $\Perspace$ 
unless the plane $V$ is rational,
that is, satisfies $\dim_\Q \left(V\cap H^2(M, \Q)\right)= 2$.

\hfill

{\bf Proof:}
Let $\Gamma\subset G$ be the monodromy group of $M$,
that is, the image of the mapping class group in $G$,
where $G=SO^+(H^2(M, \R), q)$. It is an arithmetic lattice in $G$,
as shown in \ref{_mapping_class_Theorem_}.
Since $I$ is non-ergodic, the closure $\overline{\Gamma\cdot I}$ 
of $\Gamma\cdot I$ is strictly smaller than $\Perspace$. 
By Ratner's theorem, there exists a subgroup $S\subsetneq G$ containing $H$
such that $\overline{\Gamma\cdot I}= S \cdot I$,
and $S\cap \Gamma$ is a lattice in $S$.
\ref{_subgro_SO(3,k)_Theorem_} implies that $S=H$.
Since $S\cap \Gamma$ is a lattice, this set is Zariski dense in $S$.
By \ref{_closure_arithm_Ratner_Theorem_}, $S=H$ is a rational subgroup of $G$.
Conversely, if $H$ is rational, its image is closed
in $G/\Gamma$ as follows from \ref{_closure_arithm_Ratner_Theorem_}.
\endproof

\hfill

\theorem\label{_rati_peri_tori_Theorem_}
Let $M$ be a compact complex torus of dimension $n\geq 2$,
$\Teich$ its Teichm\"uller space, $\Teich=SL(2n, \R)/SL(n,\C)$, and 
$I\in \Teich$ be a point associated to
a complex structure $I\in \End(\R^{2n})$. Then the
point $I$ is non-ergodic if and only if $H^{1,1}(M,\R)\subset H^2(M,\R)$
is a rational subspace.

\hfill

{\bf Proof:} Let $\Gamma=SL(2n,\Z)$ be the mapping
class group of $M$. A point $I$ is non-ergodic if its 
$\Gamma$-orbit in $\Teich=SL(2n, \R)/SL(n,\C)$
is not dense. By Ratner orbit closure theorem,
$\overline{\Gamma \cdot I}=S\cdot I$, where $S\supset SL(n, \C)$
is a connected Lie subgroup of $SL(2n, \R)$. Since an
intermediate subgroup $SL(2n, \R)\supset S\supset SL(n, \C)$
is equal to $SL(2n, \R)$ or $SL(n, \C)$, the point $I$
is non-ergodic if and only if $S=SL(n, \C)$ and the orbit
$\Gamma \cdot I$ is closed. By \ref{_closure_arithm_Ratner_Theorem_},
this happens if and only if the stabilizer 
$\St(I)\cong SL(n,\C)$ of $I$ is 
a rational subgroup of $SL(2n, \R)$. 
The centralizer $Z(\St(I))$ is a group 
$R_I\cong U(1)=\cos t+\sin t \cdot I$,
and $Z(Z(\St(I))=\St(I)$, hence rationality of $\St(I)$
is equivalent to rationality of $R_I$.

However, the space $H^2(M)^{R_I}$ of $R_I$-invariants is $H^{1,1}(M)$,
and, conversely, $R_I$ is a subgroup of $SL(2n, \R)=SL(H^1(M,\R))$ 
acting trivially on $H^{1,1}(M)$. Therefore, $R_I$ is rational
if and only if $H^{1,1}(M)\subset H^2(M,\R)$ is rational.
\endproof

\hfill

We have just proven density of certain orbits of $\Gamma$ in the
period space, but for geometric applications, one would need 
density of orbits in the Teichm\"uller space. This is already
true for a torus, because for the torus the period space coincides
with the Teichm\"uller space. For a hyperk\"ahler manifold
with rational curves, a similar result can be obtained directly.

\hfill

\corollary\label{_ergo_no_rat_curves_Corollary_}
Let $(M,I)$ be a hyperk\"ahler manifold with Picard of 
non-maximal rank. Assume that $(M,I)$ contains no rational
curves. Then $I$ is an ergodic complex structure,

\hfill

{\bf Proof:} Let $\Teich_0\subset \Teich$ be the set of all
Hausdorff points in $\Teich$. By \ref{_Huy_non_separa_Theorem_},
$\Teich_0$ is the set of all complex structures on $M$ admitting
non non-trivial birational models. However, any birational
map between complex manifolds with trivial canonical bundle
must blow down some subvariety, hence such map do not exist
when one has no rational curves. Therefore, $I\in \Teich_0$.
Now, the period map restricted to $\Teich_0$ is a homeomorphism,
and $\Gamma\cdot \Per(I)$ is dense in $\Perspace$ by 
\ref{_perspace_dense_Theorem_}. Therefore, $\Gamma\cdot I$
is dense in an appropriate connected 
component of $\Teich_0$, but $\Teich_0$ is dense in
$\Teich$ by \ref{_fiber_Per_over_divi_Remark_}. \endproof

\hfill

\ref{_ergo_no_rat_curves_Corollary_}
is already sufficient for many applications
dealing with hyperbolicity; indeed, to prove
that a manifold is non-hyperbolic, it suffices to show that
it contain rational curves. However, for many applications a 
full strength ergodicity result is required.

\hfill

\theorem\label{_ergo_complex_main_Theorem_}
Let $M$ be a 
hyperk\"ahler manifold, and $I$ a complex structure
of non-maximal Picard rank. Then $I$ is ergodic.

\hfill

{\bf Proof:} See Subsection \ref{_Non-Hausdorff_dense_Subsection_}. \endproof

\hfill

\corollary \label{_maxima_rank_non-ergo_Corollary_}
Let $M$ be a hyperk\"ahler manifold 
or a complex torus of complex dimension $\geq 2$. Then $M$ is non-ergodic
if and only if its Neron-Severi lattice has maximal possible
rank. In particular, there is only a countably
many non-ergodic complex structures.

\hfill

{\bf Proof:} By definition, the Neron-Severi lattice is a lattice
of integer (1,1)-classes in $H^2(M)$. It is easy to see that
it has maximal possible rank if and only if $\Per(I)$ rational
(for hyperk\"ahler manifolds). For complex tori,
the argument is given in the proof of \ref{_rati_peri_tori_Theorem_}.
The countability of the set of such
complex structures in also well known and easy to check.
\endproof

\subsection{Density of non-Hausdorff orbits}
\label{_Non-Hausdorff_dense_Subsection_}

Fix a connected component of a Teichm\"uller space
of hyperk\"ahler manifold. Abusing the notation,
we denote it $\Teich$, and denote the subgroup of the
mappoing class group fixing $\Teich$ by $\Gamma$.

Let $[I]\in \Perspace$ be a point in the period space
of $M$. The Hodge decomposition on $H^2(M)$
is determined by the periods, and we denote the corresponding (1,1)-space
by $H^{1,1}([I])$.  {\bf The positive cone}
$\Pos([I])$ is the set of all real (1,1)-classes $v\in H^{1,1}([I])$
satisfying $q(v,v)>0$. A subset $K\subset \Pos([I])$
is called {\bf K\"ahler chamber} if it is a K\"ahler cone
for some $I\in \Teich$ satisfying $\Per(I)=[I]$.
We have already used the following 
result, which is due to Eyal Markman.

\hfill

\proposition\label{_chambers_periods_teich_Proposition_}
Different K\"ahler chambers of $[I]$ do not intersect, and
$\Pos([I])$ is a closure of their union. Moreover, there is
a bijective correspondence between points of
$\Per^{-1}([I])$ in one Teichm\"uller component
and the set of K\"ahler chambers of $[I]$.

{\bf Proof:} \cite[Proposition 5.14]{_Markman:survey_}. \endproof

\hfill

Consider the set $\Hyp$ of pairs $I\in \Teich, \omega\in \Kah(M,I)$,
where $\Kah(M,I)$ denotes the K\"ahler cone. One should think of 
$\Hyp$ as of the Teichm\"uller space of all hyperk\"ahler metric
on a holomorphically symplectic manifold.  Let $F$ be the set of
all pairs $[I]\in \Perspace$, $\omega\in \Pos([I])$. 
Consider the period map $\Per_h:\; \Hyp\arrow F$ mapping
$(I,\omega)$ to $(\Per(I),\omega)$. By \ref{_chambers_periods_teich_Proposition_},
$\Per_h$ is injective with dense image. 

To prove that $\Gamma\cdot I$ is dense in $\Teich$ is the same as to show
that \[\Gamma\cdot(I,\Kah(M,I))\] is dense in $\Hyp\subset F$
(\ref{_chambers_periods_teich_Proposition_}). We consider $F$
as a homogeneous space of an appropriate Lie group. To show that 
$\Gamma\cdot(I,\Kah(M,I))$ is dense in $F$, we show that
$\Kah(M,I)$ contains an orbit of its Lie subgroup and 
apply Ratner theorem to this homogeneous space.

\hfill

Our arguments are based on the following lemma.

\hfill

\lemma\label{_Kahler_cone_contains_pos_planes_Lemma_}
Let $(M,I)$ be a hyperk\"ahler manifold,
$\omega\in \Kah(M,I)$ a K\"ahler class, 
$H^{1,1}_I(M,\Q)$ the space of rational
$(1,1)$-classes, and $l\in H^{1,1}_I(M,\Q)^\bot$
a (1,1)-class orthogonal to $H^{1,1}_I(M,\Q)$.
Then $V\cap \Pos(M,I)\subset \Kah(M,I)$,
where $V=\langle \omega, l\rangle$ is a 2-dimensional
space generated by $l$ and $\omega$.

\hfill

{\bf Proof:} As follows from \cite{_Huybrechts:cone_} 
and \cite{_Boucksom-cone_}  (see
\cite[Theorem 1.19]{_AV:MBM_} for a precise statement), 
$\Kah(M,I)$ is a subset of positive cone given by a
set of linear inequalities
\[
\Kah(M,I)=\{\omega\in \Pos(M,I)\ \ |\ \ q(\omega, l_i)> 0 \}
\]
where $l_i$ is a countable set of rational $(1,1)$-classes.
This means that for any $\omega\in \Kah(M,I)$ and any
 $l\in H^{1,1}_I(M,\Q)^\bot$, the sum $l+\omega$
also belongs to $\Kah(M,I)$, as long as
it has positive square. 
\endproof

\hfill

As we have already observed, to prove 
\ref{_ergo_complex_main_Theorem_}
it would suffice to show that $\Gamma\cdot (I,\Kah(M,I))$
is dense in $F$. Consider the set $F_1$
of all $([I],\eta)\in F$ such that $q(\eta, \eta)=1$.
Clearly, $F_1= \frac{SO(3, b_2-3)}{SO(2)\times SO(b_2-3)}$.
Indeed, $F_1$ is identified with the set of pairs
\[ 
\{ W\in \Gr_{++}(H^2(M,\R)), \omega\in W^\bot\ \ |\ \  q(\omega, \omega)>0\}.
\]

By \ref{_Kahler_cone_contains_pos_planes_Lemma_},
for any $\omega\in \Kah(M,I)$ and any 
$l\in H^{1,1}_I(M,\Q)^\bot$, the whole set
$\Pos(M,I)\cap \langle \omega, l\rangle$ belongs to 
$\Kah(M,I)$. Choose $l$ in such a way that 
$q(l,l)<0$; since $\Pic(M)$ is not of maximal rank,
this is always possible.
Consider the group $H_0$ of oriented isometries of
$V:=\langle \omega, l\rangle$; we extend its action to $H^2(M,\R)$
by requiring $H_0$ to act trivially on $V^\bot$.
By \ref{_Kahler_cone_contains_pos_planes_Lemma_}, 
$H_0$ preserves $\Kah(M,I)$. To prove density
of $\Gamma\cdot(I,\Kah(M,I))$ in $F_1$ it would suffice
to show that a $\Gamma$-orbit of the set $(I,H_0\cdot \omega)$
is dense in $F_1$. This is the same as to show that
a $\Gamma$-orbit of an appropriate point in
\[\frac{SO(3, b_2-3)}{H_0\cdot (SO(2)\times SO(b_2-3))} = \frac{SO(3, b_2-3)}{SO(2)\times SO(1,b_2-3)}
\]
is dense.

We have arrived at the situation described in
\ref{_subgro_SO(3,k)_Theorem_} and \ref{_perspace_dense_Theorem_}.
Here it was shown that any orbit of $\Gamma$ is either closed or dense.
For this orbit to be closed, the stabilizer
of a pair 
\[ [I]\in \Gr_{++}(H^2(M,\R)),\ \  V=\langle \omega, l\rangle
\]
has to be rational; since $[I]$ is irrational, this is
impossible. We proved \ref{_ergo_complex_main_Theorem_}.
\endproof

\subsection{Maximal subgroups of Lie groups}
\label{_lie_subgrou_Subsection_}

\theorem\label{_SO(p,q)_subgroup_Theorem_}
Let $(V,q)$ be a real vector space equipped with a non-degenerate
quadratic form, $G=SO^+(V)$ the connected component of the group
of isometries of $V$, and $W\subset V$ a subspace with $q\restrict W$
non-degenerate. Consider the stbgroup $H\subset G$
consisting of all isometries preserving $W\subset V$. 
Then any closed connected
Lie subgroup $S\subsetneq G$ containing $H$ coincides
with $H$.

\hfill

{\bf Proof:} 
Let ${\goth h}$, ${\goth g}$, ${\goth s}$ 
be the Lie algebras of $H$, $G$, $S$.
Then ${\goth h}=\goth{so}(W)\oplus \goth{so}(W^\bot)$.
The quotient ${\goth g}/{\goth h}$ is identified with
$\Hom(W,W^\bot)$, hence it is an 
irreducible representations of ${\goth h}$. Since 
${\goth s}/{\goth h}$ is a proper ${\goth h}$-subrepresentation of 
${\goth g}/{\goth h}$, it is equal to 0. \endproof

\hfill

\remark 
The proof of \ref{_SO(p,q)_subgroup_Theorem_}
is intuitively very clear: any isometry fixing 
$W$ is contained in $H$; if we add an isometry
which moves $W$, the resulting group should
contain all isometries. A similar argument works
for a pair $SL(n,\C)\subset SL(2n,\R)$, if we think
of $SL(n,\C)$ as of a group fixing a 
subspace of $(1,0)$-vectors in the complexification
of $\C^n$.

\hfill

\theorem\label{_complex_automo_maxi_Theorem_}
Let $W$ be a complex vector space, 
$G=SL(W,\R)$ the group of its real volume-preserving
automorphisms, and $H\cong SL(W, \C) \subset G$
the group of complex volume-preserving
automorphisms of $W$.
Then any closed connected
Lie subgroup $S\subsetneq G$ containing $H$ coincides
with $H$.

\hfill

{\bf Proof:} 
Let ${\goth h}_\C$, ${\goth g}_\C$, ${\goth s}_\C$ 
be the complexified Lie algebras of $H$, $G$, $S$.
Consider the space $W_\C:=W\otimes_\R \C$,
and let $W_\C:= W^{1,0}\oplus W^{0,1}$ be its
Hodge decomposition. Then 
\[ {\goth h}_\C={\goth {sl}}(W^{1,0})\oplus {\goth {sl}}(W^{0,1}),
\]
and 
\[ {\goth g}_\C/{\goth h}_\C= \Hom(W^{1,0},W^{0,1})\oplus 
\Hom(W^{0,1},W^{1,0}).
\] Both components
$\Hom(W^{1,0},W^{0,1})$ and $\Hom(W^{0,1},W^{1,0})$
are irreducible representations of ${\goth h}_\C$.
Since ${\goth s}_\C/{\goth h}_\C$ is a proper
${\goth h}_\C$-subrepresentation of ${\goth g}_\C/{\goth h}_\C$,
it is equal to $\Hom(W^{1,0},W^{0,1})$ or $\Hom(W^{0,1},W^{1,0})$
or 0. However, ${\goth s}$ is real, hence
${\goth s}_\C/{\goth h}_\C$ is fixed by the
anticomplex involution exchanging $W^{1,0}$ and $W^{0,1}$.
Therefore, the components $\Hom(W^{1,0},W^{0,1})$ 
and $\Hom(W^{0,1},W^{1,0})$ can be contained in 
${\goth s}_\C/{\goth h}_\C$ only together.
Since ${\goth s}_\C\subset {\goth g}_\C$ is a proper
subalgebra, ${\goth s}_\C/{\goth h}_\C$ must be empty.
\endproof


\section{Twistor spaces and Kobayashi pseudometric}
\label{_Kobayashi_Section_}

\subsection{Kobayashi pseudometric and Brody lemma}

This subsection is a brief introduction to a subject.
For more details, please see \cite{_Lang:hyperbolic_}, 
\cite{_Voisin:kobayashi_} 
and \cite{_Demailly:kobayashi_}.

\hfill

\definition
Let $M$ be a complex manifold, $x,y\in M$
points, and $d_P$ the Poincare metric on the 
unit disk $\Delta\subset \C$. Define
\[ 
  \tilde d(x,y):=
\sup\limits_{f:\; \Delta \arrow M}d_P(f^{-1}(x),f^{-1}(y))
\]
where the supremum is taken over all 
holomorphic maps $f:\; \Delta \arrow M$
from the disk $\Delta$ to $M$ 
such that $f(\Delta)\supset \{x,y\}$.
The maximal pseudo-metric $d$ satisfying
$d(x,y) \leq \tilde d(x,y)$
is called {\bf the Kobayashi pseudometric}.
The manifold $M$ is called {\bf Kobayashi hyperbolic}
if the Kobayashi pseudometric is non-degenerate
(\cite{_Kobayashi:1970_}).

\hfill

For a compact manifold, hyperbolocity can be interpreted
as non-\-exis\-tence of entire curves.

\hfill

\definition
{\bf An entire curve} in a complex manifold $M$
is an image of a non-constant holomorphic map $\C \arrow M$.

\hfill

The following two theorems 
are fundamental in hyperbolic geometry;
for details and the proofs, see
again \cite{_Lang:hyperbolic_},
\cite{_Voisin:kobayashi_} and \cite{_Demailly:kobayashi_}.
They follow from a remarkable result on convergence
of disks and entire curves on complex manifolds, 
called Brody's lemma (\cite{_Brody:hyperboloc_}).

\hfill

\theorem
Let $M$ be a compact complex manifold.
Then $M$ contains an entire curve if and only if
it is not Kobayashi hyperbolic. 
\endproof

\hfill

\theorem\label{_limi_non-hy_Theorem_}
Let $I_i$ be a sequence of non-hyperbolic 
complex structures on 
a compact manifold $M$, and $I$ its limit. Then
$(M,I)$ is also non-hyperbolic.
\endproof

\hfill

The main result of this section is the following theorem.

\hfill

\theorem\label{_non_hype_Theorem_}
Any compact hyperk\"ahler manifold 
$M$ with $b_2(M)> 3$ is non-hyperbolic.

\hfill

{\bf Proof:} See Subsection 
\ref{_twi_campa_Subsection_}. \endproof

\hfill

\remark 
For all known examples of hyperk\"ahler manifolds,
this theorem is already known, due to Kamenova
and Verbitsky (\cite{_Kamenova_V:fibrations_}).

\hfill

\remark 
To prove \ref{_non_hype_Theorem_},
it would suffice to show that
there exists an ergodic complex structure $I$
which is non-hyperbolic. Indeed, in this case
the orbit of $I$ is dense. This implies that any complex
structure can be obtained as a limit of non-hyperbolic ones
(\ref{_ergo_limits_Remark_}).

\subsection{Twistor spaces and Campana theorem}
\label{_twi_campa_Subsection_}

\definition 
Let $I,J,K,g$ be a hyperk\"ahler structure on a manifold $M$.
{\bf Induced complex structures} on $M$ are 
complex structures of form 
$S^2 \cong \{ L:= aI + bJ +c K, \ \ \ a^2+b^2+c^2=1.\}$

\hfill

\definition\label{_twistor_Definition_}
A {\bf twistor space} $\Tw(M)$ of a hyperk\"ahler manifold
is a complex manifold obtained by gluing induced complex structures into
a holomorphic family over $\C P^1$. More formally:

{\em Let $\Tw(M) := M \times S^2$. 
Consider the complex structure $I_m:T_mM \to T_mM$ 
on $M$ induced by $J \in S^2 \subset {\Bbb H}$. Let $I_J$
denote the complex structure on $S^2 = \C P^1$.
The operator $I_{\Tw} = I_m \oplus I_J:T_x\Tw(M) \to T_x\Tw(M)$ 
satisfies $I_{\Tw} ^2 = -\Id$. It defines 
an almost complex structure on $\Tw(M)$. This almost 
complex structure is known to be integrable 
(\cite{_Obata_} and \cite{_Salamon_}; 
see \cite{_Kaledin:twistor_} for a modern proof).}

\hfill

Rational curves on twistor spaces were studied 
by F. Campana in a series of papers (\cite{_Campana:rational_lines_},
\cite{_Campana:twistors_}); among the results
of this study, Campana proved the following theorem.

\hfill

\theorem \label{_Campana_non_hy_Theorem_}
(\cite{_Campana:twistor_hy_})
Let $M$ be a hyperk\"ahler manifold, equipped with a hyperk\"ahler 
structure, and 
$\Tw(M)\stackrel \pi \arrow \C P^1$ its twistor
space. Then there exists an entire curve in some
fiber of $\pi$. \endproof

\hfill

\claim\label{_twi_ergo_Claim_}
Let $M$ be a hyperk\"ahler manifold, $b_2(M)\geq 4$.
Then there exists a twistor family on $M$ which has only ergodic
fibers. 

\hfill

{\bf Proof:} 
By \ref{_maxima_rank_non-ergo_Corollary_}, 
there are only countably many complex
structures which are not ergodic. The space ${\cal T}$ of all
twistor families is identified with the set of hyperk\"ahler metrics
up to a constant multiplier. Therefore, it has real dimension 
$\frac{b_2(b_2-1)(b_2-2)}6$, as follows from Bogomolov's
local Torelli theorem and the Calabi-Yau theorem
(\ref{_Calabi_Yau_Theorem_}).
The space of twistor families passing through a
given complex structure is parametrized by the
projectivization of a K\"ahler cone, hence its real
dimension is $b_2-3$. There is a countable number
of non-ergodic complex structures, hence the set
${\cal T}_0$ of twistor families passing through non-ergodic
complex structures is a union of countably many
$b_2-3$-dimensional families. For $b_2>3$, one has
$\dim_\R {\cal T}> b_2 -3$,  hence 
${\cal T}_0$ has measure 0 in
${\cal T}$. \endproof

\hfill

Non-hyperbolicity of a hyperk\"ahler manifold
follows from this claim immediately.
Indeed, let $\pi:\; {\cal S}\arrow \C P^1$ be a twistor family
with all fibers ergodic.
\ref{_Campana_non_hy_Theorem_}
implies that at least one fiber of $\pi$ is non-hyperbolic.
Denote this fiber by $M$. Since $M$ is ergodic,
there is a dense family of manifolds biholomorphic
to $M$ in the Teichm\"uller space $\Teich$. Since non-hyperbolic 
complex structures are closed in $\Teich$ (\ref{_limi_non-hy_Theorem_}), 
this implies that all 
points in $\Teich$ correspond to non-hyperbolic
complex structures.

\hfill

{\bf Acknowledgements:} I am grateful to 
Frederic Campana, Ljudmila Kamenova, Eyal Markman,
Ossip Shvartsman and Sasha Anan$'$in for interesting discussions 
of the subject. Dima Panov was the first to 
point out the pathological behaviour of the
mapping group action on the Teichm\"uller space.
This paper would not appear
without Dmitry Kleinbock, who has explained to 
me how to apply the ergodic theory to the subject,
and gave an ample reference to the Moore's and 
Ratner's theorems. I am grateful to Konstantin
Tolmachev for interesting talks on ergodic theory
and its applications. Many thanks to Eyal Markman
for finding a gap in one of the last versions of
this paper (fixed in Subsection \ref{_Non-Hausdorff_dense_Subsection_}).

\hfill

{\small
\noindent {\sc Misha Verbitsky\\
{\sc Laboratory of Algebraic Geometry,\\
National Research University HSE,\\
Faculty of Mathematics, 7 Vavilova Str. Moscow, Russia,}\\
\tt  verbit@mccme.ru}, also: \\
{\sc Kavli IPMU (WPI), the University of Tokyo}
}


{\scriptsize
\begin{thebibliography}{AV1}


\bibitem[AmV]{_AV:MBM_} 
Amerik, E., Verbitsky, M. {\em Rational curves on 
hyperk\"ahler manifolds}, arXiv:1401.0479

\bibitem[AnV]{_Ananin_Verbitsky_}
Sasha Anan'in, Misha Verbitsky,
{\em Any component of moduli of polarized hyperk\"ahler
manifolds is dense in its deformation space},
arXiv:1008.2480, 17 pages, 4 figures.


\bibitem[Bea]{_Beauville_} 
 Beauville, A. {\em 
Varietes K\"ahleriennes dont la premi\`ere classe de Chern est
nulle.}  J. Diff. Geom. {\bf 18} (1983) 755-782.

\bibitem[Bes]{_Besse:Einst_Manifo_} 
Besse, 
A., {\em Einstein Manifolds}, Springer-Verlag, New York (1987)

\bibitem[Bo1]{_Bogomolov:decompo_}  
Bogomolov, F. A., {\em On the decomposition of 
K\"ahler manifolds with trivial canonical class}, Math. USSR-Sb.
{\bf 22} (1974) 580 - 583.

\bibitem[Bo2]{_Bogomolov:defo_} 
F. Bogomolov, {\em Hamiltonian K\"ahler
manifolds}, Sov. Math. Dokl. {\bf 19} (1978) 1462 - 1465.


\bibitem[Bou1]{_Boucksom-cone_} 
Boucksom, S., {\em Le c\^one k\"ahl\'erien d'une vari\'et\'e 
hyperk\"ahl\'erienne}, C. R. Acad. Sci. Paris 
Ser. I Math. 333 (2001), no. 10, 935--938.


\bibitem[Br]{_Brody:hyperboloc_}
 Brody, R., {\em Compact manifolds and hyperbolicity,}
Trans. Amer. Math. Soc. 235 (1978), 213-219.



\bibitem[Cam1]{_Campana:rational_lines_}
Campana, F., {\em
Alg\'ebricit\'e et compacit\'e dans les espaces de cycles}, 
C. R. Acad. Sci. Paris S\'er. A-B 289 (1979), no. 3, A195eA197. 

\bibitem[Cam2]{_Campana:twistors_}
Campana, F.,
{\em On twistor spaces of the class C,}
J. Differential Geom. 33 (1991), no. 2, 541-549. 


\bibitem[Cam3]{_Campana:twistor_hy_}
 F. Campana,
{\em An application of twistor theory to 
the nonhyperbolicity of certain compact symplecti
c K\"ahler manifolds} 
J. Reine Angew. Math., 425:1-7, 1992.



\bibitem[Cat]{_Catanese:moduli_}
F. Catanese, 
{\em A Superficial Working Guide to Deformations and Moduli},
arXiv:1106.1368, 
Advanced Lectures in Mathematics, Volume XXVI
Handbook of Moduli, Volume III, page 161-216 
(International Press).


\bibitem[D]{_Demailly:kobayashi_}
Demailly, Jean-Pierre,
{\em Hyperbolic algebraic varieties and holomorphic differential equations,}
 Talk given at the VIASM annual meeting, 25-26 August
 2012, long expository version (135 pages),
\url{http://www-fourier.ujf-grenoble.fr/~demailly/manuscripts/viasm2012expanded.pdf}.


\bibitem[DDM]{_DD_Mourougane_}
Guillaume Deschamps, No\"el Le Du, Christophe Mourougane,
{\em Hessian of the metric form on twistor spaces},
arXiv:1202.0183.


\bibitem[Gr]{_Gromov:Riemannian_} 
Gromov, Misha, {\em  Metric structures for Riemannian and
non-Riemannian spaces}, Based on the 1981 French
original. With appendices by M. Katz, P. Pansu and
S. Semmes. Translated from the French by Sean Michael
Bates. Progress in Mathematics, 152. Birkh\"auser Boston,
Inc., Boston, MA, 1999. xx+585 pp.


\bibitem[GTZ1]{_Gross_Tosatti_Zhang:2011_}
M. Gross, V. Tosatti, Y. Zhang, {\em Collapsing 
of Abelian Fibered Calabi-Yau Manifolds,}
Duke Math. J. 162 (2013), no. 3, 517-551.


\bibitem[GTZ2]{_Gross_Tosatti_Zhang:2013_}
Mark Gross, Valentino Tosatti, Yuguang Zhang,
{\em Gromov-Hausdorff collapsing of Calabi-Yau manifolds},
arXiv:1304.1820.


\bibitem[F]{_Fujiki:HK_}  
Fujiki, A. {\em On the de Rham Cohomology Group of a Compact 
K\"ahler Symplectic Manifold}, Adv. Stud.
Pure Math. 10 (1987) 105-165.


\bibitem[HKLR]{_HKLR_}  
N. J. Hitchin, A. Karlhede, 
U. Lindstr\"om, M. Ro\v cek, 
{\em Hyperk\"ahler metrics and supersymmetry}, 
Comm. Math. Phys. {\bf 108} (1987), 535-589.



\bibitem[H1]{_Huybrechts:basic_} 
Huybrechts, D., 
{\em Compact hyperk\"ahler manifolds: Basic
results}, Invent. Math. 135 (1999), 63-113, 
alg-geom/9705025.


\bibitem[H2]{_Huybrechts:erratum_} 
Huybrechts, D., 
{\em Erratum to the paper: Compact hyperk\"ahler
  manifolds: basic results},
 Invent. math. 152  (2003), 209-212, math.AG/0106014.


\bibitem[Hu3]{_Huybrechts:lec_}
 Huybrechts, Daniel, 
{\em Compact hyperk\"ahler manifolds, Calabi-Yau 
manifolds and related geometries,}
Universitext, Springer-Verlag, Berlin, 2003, 
Lectures from the Summer School held in
Nordfjordeid, June 2001, pp. 161-225.



\bibitem[H4]{_Huybrechts:cone_}
Huybrechts, D.,
{\em The K\"ahler cone of a compact hyperk\"ahler manifold},
Math. Ann. 326 (2003), no. 3, 499--513, arXiv:math/9909109.


\bibitem[Ka]{_Kaledin:twistor_}
D. Kaledin,
{\em Integrability of the twistor space for a hypercomplex manifold},
Selecta Math. (N.S.) {\bf 4} (1998) 271-278.



\bibitem[KaV]{_NHYM_} 
Kaledin, D., Verbitsky, M.,
{\it Non-Hermitian 
Yang-Mills connections}, Selecta Math. (N.S.) {\bf 4}
(1998), no. 2, 279--320.


\bibitem[KV]{_Kamenova_V:fibrations_}
Ljudmila Kamenova, Misha Verbitsky 
{\em Families of Lagrangian fibrations on hyperk\"ahler
  manifolds}, arXiv:1208.4626, 13 pages.



\bibitem[KSS]{_Kleinbock_etc:Handbook_}
Kleinbock, Dmitry; Shah, Nimish; Starkov, Alexander,
{\em Dynamics of subgroup actions on homogeneous spaces of
  Lie groups and applications to number theory}, Handbook
of dynamical systems, Vol. 1A, 813-930, North-Holland,
Amsterdam, 2002. 

\bibitem[Ko]{_Kobayashi:1970_}
Kobayashi, Shoshichi,
{\em Hyperbolic manifolds and holomorphic mappings,}
 Pure and Applied Mathematics 2,  (1970), New York: Marcel Dekker Inc.

\bibitem[KS]{_Kontsevich_Soibelman:torus_}
 M. Kontsevich, Y. Soibelman, 
{\em Homological mirror symmetry and torus fibrations,} in
 Symplectic geometry and mirror symmetry, 203-263, World
 Sci. Publishing 2001.



\bibitem[L]{_Lang:hyperbolic_}
S. Lang, {\em 
Introduction to complex hyperbolic spaces}, Springer, New York, 1987.


\bibitem[Ma1]{_Markman:constra_}
Markman, E. {\em
Integral constraints on the monodromy group of
    the hyperkahler resolution of a symmetric product of a
    K3 surface,} International Journal of Mathematics
Vol. 21, No. 2 (2010) 169-223, arXiv:math/0601304.


\bibitem[Ma2]{_Markman:survey_}
Markman, E. {\em
A survey of Torelli and monodromy results for 
holomorphic-symplectic varieties}, 
Proceedings of the conference "Complex and Differential 
Geometry'', Springer Proceedings in Mathematics, 2011, Volume 8, 257--322,
arXiv:math/0601304.

\bibitem[MM]{_Markman_Mehrota_}
Eyal Markman, Sukhendu Mehrotra,
{\em Hilbert schemes of K3 surfaces are dense in moduli},
11 pages, arXiv:1201.0031.


\bibitem[Mor]{_Morris:Ratner_}
Morris, Dave Witte, 
{\em Ratner's Theorems on Unipotent Flows,} 
Chicago Lectures in Mathematics, University 
of Chicago Press, 2005.


\bibitem[Mo]{_Moore:ergodi_}
Calvin C. Moore,
{\em Ergodicity of Flows on Homogeneous Spaces},
American Journal of Mathematics
Vol. 88, No. 1 (Jan., 1966), pp. 154-178


\bibitem[Ob]{_Obata_}
Obata, M., {\em Affine connections 
on manifolds with almost complex, quaternionic or Hermi
tian structure}, Jap.
J. Math., 26 (1955), 43-79.

\bibitem[PS]{_Piatetski_Shapiro_Shafarevich_}
Piatecki-Shapiro, I.I.; Shafarevich I.R.,
{\em Torelli's
theorem for algebraic surfaces of type
K3}, Izv. Akad. Nauk SSSR Ser. Mat. (1971) 35: 530-572.



\bibitem[Sal]{_Salamon_} 
S. Salamon, 
{\em Quaternionic K\"ahler manifolds,}
Inv. Math. {\bf 67} (1982) 143--171.


\bibitem[Sh]{_Shah:uniformly_}
N. A. Shah,
{\em Uniformly distributed orbits of certain flows on homogeneous spaces,} 
Math. Ann. 289 (2) (1991), 315-33.



\bibitem[Su]{_Sullivan:infinite_}
Sullivan, D.,
{\em Infinitesimal computations in topology}, Publications
Math\-\'ema\-tiques de l'IH\'ES, 47 (1977), p. 269-331

\bibitem[V1]{_Verbitsky:coho_announce_} 
Verbitsky, M.,
{\it Cohomology of compact hyperk\"ahler manifolds
and its applications,} GAFA vol. 6 (4) pp. 601-612 (1996).



\bibitem[V2]{_V:Torelli_}
Verbitsky, M.,
{\em A global Torelli theorem for hyperk\"ahler manifolds,}
 Duke Math. J. \textbf{162} (2013), 2929-2986.


\bibitem[V3]{_Verbitsky:special_twistor_}
Misha Verbitsky,
{\em Rational curves and special metrics on twistor spaces},
 arXiv:1210.6725, 12 pages.


\bibitem[V4]{_Verbitsky:hypercomple_}
Verbitsky M., {\em Hypercomplex Varieties}, 
 Comm. Anal. Geom. {\bf 7} 
(1999), no. 2, 355--396.


\bibitem[Vi]{_Viehweg:moduli_}
Viehweg, E.,
{\em Quasi-projective Moduli for Polarized Manifolds,}
Springer-Verlag, Berlin, Heidelberg, New York, 1995,
Ergebnisse der Mathematik und ihrer Grenzgebiete, 3.~Folge, Band 30,
also available at
{\tt http://www.uni-due.de/$\widetilde{\phantom{a}}$mat903/books.html}

\bibitem[Vo]{_Voisin:kobayashi_}
Claire Voisin,
{\em On some problems of Kobayashi and Lang; algebraic
approaches,}  Current Developments in Mathematics 2003,
no. 1 (2003), 53-125.


\bibitem[VGS]{_Vinberg_Gorbatsevich_Shvartsman_}
Vinberg, E. B.,  Gorbatsevich, V. V.,  Shvartsman, O. V., 
{\em Discrete Subgroups of Lie Groups}, in 
``Lie Groups and Lie Algebras II'',
Springer-Verlag, 2000.


\bibitem[Y]{_Yau:Calabi-Yau_} 
Yau, S. T., {\em On the Ricci curvature of a compact K\"ahler manifold 
and the complex Monge-Amp\`ere equation I.} Comm. on Pure and Appl.
Math. 31, 339-411 (1978).



\end{thebibliography}
}
\end{document}